\documentclass[12pt]{article}

\usepackage{amsthm,amsmath,amssymb}

\usepackage{graphicx}

\usepackage[colorlinks=true,citecolor=black,linkcolor=black,urlcolor=blue]{hyperref}

\theoremstyle{plain}
\newtheorem{theorem}{Theorem}
\newtheorem{lemma}[theorem]{Lemma}
\newtheorem{corollary}[theorem]{Corollary}

\theoremstyle{definition}

\newtheorem{conjecture}[theorem]{Conjecture}

\theoremstyle{remark}

\begin{document}

\title{The extremal function for cycles of length $\ell$ mod $k$}

\author{Benny Sudakov\thanks{Research supported by SNSF grant 200021-149111 and a USA-Israel BSF grant.}\\
\small Department of Mathematics \\ [-0.8ex]
\small ETH Z\"{u}rich \\ [-0.8ex]
\small Ramistrasse 101 \\ [-0.8ex]
\small 8092 Z\"{u}rich, Switzerland. \\
\small\tt benjamin.sudakov@math.ethz.ch. \\
\and
Jacques Verstra\"ete\thanks{Research supported by the National Science Foundation
and by the Institute for Mathematical Research (FIM) of ETH Z\"urich.} \\
\small Department of Mathematics \\ [-0.8ex]
\small University of California at San Diego\\ [-0.8ex]
\small 9500 Gilman Drive, La Jolla \\ [-0.8ex]
\small California 92093-0112, USA. \\
\small\tt jverstra@math.ucsd.edu. }

\maketitle

\begin{abstract}
Burr and Erd\H{o}s conjectured that for each $k,\ell \in \mathbb Z^+$ such that $k \mathbb Z + \ell$ contains even integers, there exists $c_k(\ell)$ such that any graph of average degree at least $c_k(\ell)$ contains a cycle of length $\ell$ mod $k$. This conjecture was proved by Bollob\'{a}s, and many successive improvements of upper bounds on $c_k(\ell)$ appear in the literature. In this short note, for $1 \leq \ell \leq k$, we show that $c_k(\ell)$ is proportional to the largest average degree of a $C_{\ell}$-free graph on $k$ vertices, which determines $c_k(\ell)$ up to an absolute constant. In particular, using known results on Tur\'{a}n numbers for even cycles, we obtain $c_k(\ell) = O(\ell k^{2/\ell})$
for all even $\ell$, which is tight for $\ell \in \{4,6,10\}$. Since the complete bipartite graph $K_{\ell - 1,n - \ell + 1}$ has no cycle of length $2\ell$ mod $k$, it also shows $c_k(\ell) = \Theta(\ell)$ for $\ell = \Omega(\log k)$.
\end{abstract}

\section{Introduction}

More than forty years ago, Burr and Erd\H{o}s conjectured that for each $k,\ell \in \mathbb Z^+$ such that $k \mathbb Z + \ell$ contains even integers,  there exists $c_k(\ell)$ such that any graph of average degree at least $c_k(\ell)$ contains a cycle of length $\ell$ mod $k$. Bollob\'{a}s~\cite{Bollobas} was the first to show that $c_k(\ell)$ exists, showing $c_k(\ell) \leq \frac{1}{2k}[(k + 1)^{k} - 1]$ for all $\ell$. This upper bound was reduced by a number of authors~\cite{BondyVince,Chen,Dean1,Dean2,Diwan,Fan,GyarfasLehel,KostochkaVerstraete,Thomassen}. The first linear bound $c_k(\ell) \leq 8(k - 1)$ was given by the second author~\cite{Verstraete}. Thomassen~\cite{Thomassen} conjectured
that every graph of minimum degree at least $k + 1$ contains cycles of all possible even lengths mod $k$. Since every graph with average degree a least $2k$ contains a subgraph of minimum degree at least $k + 1$, this conjecture implies that $c_k(\ell) \leq 2k$ when $k$ is even.
Recently, Liu and Ma~\cite{LiuMa} showed for $k$ even that a graph of minimum degree at least $k + 1$ contains cycles of all possible even lengths mod $k$, answering a conjecture of Thomassen~\cite{Thomassen} for even values of $k$. Since any graph
whose blocks are cliques of order $k + 1$ contains no cycle of length 2 modulo $k$, this is best possible.

\subsection{Main Theorem}

The aim of this note is to determine $c_k(\ell)$ up to a constant factor. In particular we show perhaps surprisingly that $c_k(\ell)$ is  sublinear in $k$ if $\ell$ is even and $\ell = o(k)$ and in fact linear in $\ell$ if $\ell = \Omega(\log k)$ as $k \rightarrow \infty$. A number of bounds on $c_k(\ell)$ have been given in the literature~\cite{BondyVince,Dean1,Dean2,Diwan,HS} for specific values of $k$ and $\ell$. Let $d_k(\ell)$ be the
largest possible average degree of any $k$-vertex bipartite $C_{\ell}$-free graph.
Our main theorem
shows that $c_k(\ell)$ is proportional to $d_k(\ell)$.

\begin{theorem}\label{main}
For $3 \leq \ell < k$, every $C_{\ell}$-free graph of average degree at least $96 \cdot d_k(\ell)$ contains cycles of $k$ consecutive even lengths. In particular,
 \[ d_k(\ell) \leq c_k(\ell) \leq 96 \cdot d_k(\ell).\]
\end{theorem}

Since a graph on $k$ vertices with no cycle $C_{\ell}$ has no cycle of length $\ell$ mod $k$, $c_k(\ell) \geq d_{k}(\ell)$, providing the lower bound in Theorem \ref{main}. In fact one can do slightly better, by considering a large graph whose blocks are all extremal $C_{\ell}$-free graphs with $k + \ell - 1$ vertices. This gives $c_{k}(\ell) \geq d_{k + \ell - 1}(\ell)+\frac{2}{k+\ell}$. We did not attempt to optimize the
upper bound in the proof of this theorem, and the constant 96 can no doubt be improved. However, we make the following conjecture which suggests the constant in front of $d_k(\ell)$ should be $1 + o_k(1)$.

\begin{conjecture}\label{cd}
 For each even integer $\ell \geq 2$, $c_k(\ell) \sim d_k(\ell)$ as $k \rightarrow \infty$.
\end{conjecture}

\subsection{Quantitative results}

To state quantitative results, if $\ell$ is odd,
then $ \frac{k^2-1}{2k}  \leq d_k(\ell) \leq k/2$ (the average degree in the densest $k$-vertex complete bipartite graph) and so $c_k(\ell) = \Theta(k)$.
The bound $d_k(\ell) = O(\ell k^{2/\ell})$ for $\ell$ even follows from known bounds on Tur\'{a}n numbers for even cycles, the first of which were obtained by Bondy and Simonovits~\cite{BS}. The current best upper bound on $d_k(\ell)$ for fixed $\ell$ is due to Bukh and Jiang~\cite{BukhJiang}. This gives the following corollary:

\begin{corollary}
For $1 \leq \ell \leq k$, $c_k(\ell) = \Theta(k)$ if $\ell$ is odd, and $c_k(\ell) = O(\ell k^{\frac{2}{\ell}})$ if $\ell$ is even.
\end{corollary}

In particular this shows $c_k(\ell) = \Theta(\ell)$ for $\ell = \Omega(\log k)$, since a complete bipartite graph with one part of size less than $\ell/2$ contains no cycle of length $\ell$ mod $k$. Since $d_k(\ell) = \Theta(k^{2/\ell})$ when $\ell \in \{4,6,10\}$ (see~\cite{Vsurvey} and F\"{u}redi and Simonovits~\cite{FS}), we find
\[ c_k(4) = \Theta(k^{1/2}) \quad \quad c_k(6) = \Theta(k^{1/3}) \quad \quad c_k(10) = \Theta(k^{1/5})\]
so the corollary is tight up to constants when $\ell \in \{4,6,10\}$ and for $\ell = \Omega(\log k)$.
For cycles of length $4$ mod $k$ this substantially improves a result of Diwan~\cite{Diwan},
stating that a graph of minimum degree at least $k + 1$ contains a cycle of length four modulo $k$.

\subsection{Remarks on chromatic number and cycles mod $k$}

 One may attempt to extend Theorem \ref{main} to obtain cycles of any length modulo $k$ under conditions on the chromatic number of the graph. Let $\chi_k(\ell)$ denote the maximum possible chromatic number of a graph with no cycle of length $\ell$ mod $k$, and let
$n_k(\ell)$ denote the largest chromatic number in a $C_{\ell}$-free  graph on $k$ vertices. Bounds on $n_k(\ell)$ come from bounds on
cycle-complete graph Ramsey numbers. In particular, the results of Caro and Yuster~\cite{CY} and the results of the first
author~\cite{Sudakov} show that if $m = \lceil \ell/2 \rceil$, then for $\ell \geq 3$,
\[ n_k(\ell) = \Omega\Bigl(\frac{k^{1/m}}{(\log k)^{1/m}}\Bigr).\]
This is known to be tight up to a constant factor only for $\ell = 3$, for in this case Kim~\cite{Kim} constructed
$k$-vertex triangle-free graphs with chromatic number $O(\sqrt{k/\log k})$. The analog of Theorem \ref{main} for chromatic number
is stated in~\cite{KSV}, and in particular, the following result holds:

\begin{theorem}\label{main2} {\rm \cite{KSV}}
There exists $c > 0$ such that for $k \geq \ell \geq 1$ and $m = \lceil \ell/2 \rceil$,
\[ n_k(\ell) \leq \chi_k(\ell) \leq c\Bigl(\frac{k^{1/m}}{(\log k)^{1/m}}\Bigr).\]
\end{theorem}

This theorem is a special case of Theorem 4 in~\cite{KSV}, which shows that a graph of chromatic number $d$ with no cycle of length $\ell$ contains cycles of $\Omega(d^{m}\log d)$ consecutive lengths, and so in particular if $\chi_k(\ell) = d$ then $d^m \log d = O(k)$.

\medskip

Let $\delta_k(\ell)$ denote the largest minimum degree in a 2-connected non-bipartite graph with no cycle of length $\ell$ mod $k$. For $r = \lfloor \ell/2\rfloor$, it follows from results of Ma~\cite{Ma} that for each fixed $\ell \geq 3$, $\delta_k(\ell) = O(k^{1/r})$, so when $\ell$ is odd, Theorem \ref{main2} offers a stronger conclusion for graphs of large chromatic number. For instance in the case $\ell = 3$, Theorem \ref{main2} gives
\[ \chi_k(3) = \Theta\Bigl(\frac{k^{1/2}}{(\log k)^{1/2}}\Bigr),\]
whereas Ma's result gives $\delta_k(3) = \Theta(k)$. Bondy and Vince~\cite{BondyVince} showed that there exist 2-connected non-bipartite graphs of arbitrarily large minimum degree with no cycles of lengths differing by exactly 1. For 3-connected graphs, Fan~\cite{Fan} showed that every non-bipartite 3-connected graph with minimum degree at least $3k$ contains $2k$ cycles of consecutive lengths.  We may conjecture the following analog of Conjecture \ref{cd}:

\begin{conjecture}
For all $\ell \geq 3$, $\chi_k(\ell) \sim n_k(\ell)$ as $k \rightarrow \infty$.
\end{conjecture}

We remark that since there is a wide gap between upper and lower bounds on cycle-complete graph Ramsey numbers in general,
the actual asymptotic value of $n_k(\ell)$ is likely to be very difficult to determine.

\section{Preliminaries}
In this section we collect some results which we will use in our proofs.

\bigskip
\bigskip

\noindent
{\bf Breadth first search trees:} \quad If $v$ is a vertex of a connected graph $G$, then
 a {\em breadth first search tree rooted at $v$} is a spanning tree $T$ of $G$ created as follows.
 Having found a tree $T_i \subseteq G$ with vertices $v = v_0,v_1,v_2,\dots,v_i$, we pick a vertex $v_{i + 1} \in V(G) \backslash V(T_i)$
 such that the distance $d_G(v,v_{i + 1})$ is a minimum, and then let $V(T_{i + 1}) = V(T_i) \cup \{v_{i + 1}\}$ and
 select a vertex $v_j \in V(T_i)$ such that $\{v_j,v_{i + 1}\} \in E(G)$ and $j$ is a minimum,
 and set $E(T_{i + 1}) = E(T_i) \cup \{v_j,v_{i + 1}\}$. If $T$ is a breadth first search tree in $G$, rooted at $v$,
 then $d_G(v,w) = d_T(v,w)$ for all $w \in V(G)$ -- in other words $T$ preserves the distance from $v$ in $G$.
 In particular, the {\em $i$th level of $T$} is $L_i(T) = \{w \in V(G) : d_G(v,w) = i\}$. The {\em height}
 of $T$ is $\max\{d_T(v,w) : w \in V(T)\}$. Note that the edges of $G$ lie either between two consecutive levels of $T$, or inside levels of $T$.

\bigskip
\bigskip

\noindent
{\bf Theta graphs:} \quad
A {\em $\theta$-graph} consists of a cycle plus an additional edge joining two non-consecutive vertices on the cycle.
The following lemma due to Bondy and Simonovits~\cite{BS} is required for the proofs of Theorems \ref{main} and \ref{main2}.

 \begin{lemma}\label{theta1}
 Let $G$ be an $n$-vertex $\theta$-graph and let $A \cup B$ be a partition of vertices of $G$ into two non-empty subsets. Then for every $r < n$,
 $G$ contains a path of length $r$ with one end in $A$ and one end in $B$, unless
 $G$ is bipartite graph with parts $A$ and $B$.
 \end{lemma}

\noindent
 The proof of this lemma is given in~\cite{Verstraete} (see Lemma 2 in~\cite{Verstraete}). The purpose of the additional edge $e$ in a $\theta$-graph is to preclude the possibility that the vertices of $A$ occur at every $m$th vertex along the cycle $G - e$ for some $m|n$, for in that case,  there is no path of length zero mod $m$ with one end in $A$ and one end in $B$.

\bigskip
\bigskip

\noindent
{\bf Long cycles in $C_{\ell}$-free graphs:} \quad
To find $\theta$-graphs we use the following two lemmas.

 \begin{lemma}\label{theta2}
 For $3 \leq \ell \leq k$, if $F$ is a bipartite graph of average degree at least $24d_k(\ell)$ containing no $C_{\ell}$, then $F$ contains a
 $\theta$-graph with at least $2k + 2$ vertices.
 \end{lemma}

 \noindent
 {\bf Proof.}\, The proof of this lemma uses some ideas from~\cite{SV}. Since the average degree of $F$ is at least $24d_k(\ell)$, it contains a subgraph $F'$ of minimum degree at least $12d_k(\ell)$. Let $X$ be a subset of $F'$ of size $t \leq k$. We claim that
 $|\partial X| > 2|X|$, where $\partial X$ is a set of all vertices of $F' \setminus X$ which have at least one neighbor in $X$. Indeed, if this is not the case then the subset $Y=X \cup \partial X$ has size at most $3t$ and contains all the edges of $F'$ incident with $X$. Thus the average degree of the induced subgraph $F'[Y]$ is at least $12 t d_k(\ell)/(3t)=4d_k(\ell)$. If $3t\geq k$, by taking a random subset of $Y$ of size $k$ and using linearity of expectation we obtain a subgraph with average degree at least $\frac{k}{|Y|}4d_k(\ell) \geq \frac{4}{3}d_k(\ell)>d_k(\ell)$. Such graph has an $\ell$-cycle, contradiction. If $3t <k$, we can take $\lfloor k/(3t) \rfloor$ disjoint copies
 of graph $F'[Y]$ together with at most $k/2$ isolated vertices to get a $C_{\ell}$-free $k$-vertex graph with average degree
 at least $4d_k(\ell)/2>d_k(\ell)$, contradicting the definition of $d_k(\ell)$.

To complete the proof we use a variant of P\'{o}sa's rotation lemma: Lemma 2.1 in \cite{PS} states that if $|\partial X| > 2|X|$ for every subset $X \subset V(F')$ of size at most $k$, then $F'$ contains a cycle of length $3k \geq 2k+2$. Moreover one vertex of this cycle has all its neighbors on the cycle, so since the minimum degree of $F'$ is at least $12d_k(\ell)\geq 3$, this gives us a $\theta$-graph with at least $2k + 2$ vertices. \hfill $\Box$

\section{Proof of Theorem \ref{main}}

Let $G$ be a $C_{\ell}$-free graph of average degree at least $96d_k(\ell)$.
Consider a connected bipartite subgraph $G'$ of $G$ of
average degree at least $48d_k(\ell)$. Let $T$ be a breadth first search tree in $G'$. Since $G'$ is bipartite
\[ e(G') = \sum_{i\geq 1}e(L_{i-1}(T),L_i(T)).\]
On the other hand
\[ e(G')=24d_k(\ell)|V(G')| = 24d_k(\ell)\sum_{i\geq 0} |L_i(T)| \geq 12d_k(\ell)\sum_{i\geq i}(|L_{i-1}(T)|+|L_i(T)|).\]
Therefore for some $i$, the edges of $G'$ between $L_{i - 1}(T)$
and $L_i(T)$ form a bipartite graph $F$ of average degree at least $24d_{k}(\ell)$. By Lemma \ref{theta2}, $F$ contains a
$\theta$-graph $F'$ with at least $2k + 2$ vertices. Let $U = V(F') \cap L_{i - 1}(T)$ and
$W =  V(F') \cap L_i(T)$. Let $T'$ be the minimum subtree of $T$ such that $V(T') \cap U = U$.
Then the vertex $u$ of $T'$ closest to the root of $T$ has degree at least two in $T'$. This implies $T' - \{u\}$
has at least two components. Let $A$ be the set of vertices of $U$ in one of the components,
and $B = V(F') \backslash A$. Then $V(F')$ has a partition $A \cup B$, but
$A$ and $B$ do not form the bipartition of $F'$. By Lemma \ref{theta1}, for each $r \in \{1,2,\dots,k\}$, there exists a path $P$ of length $2r$ in $F'$ with one end $a \in A$ and one end $b \in B$. Since $P$ has even length,
$b \in V(T') \cap L_{i - 1}(T)$. Let $h$ denote the height of $T'$. Then $a,b \in V(T')$ are connected in $T'$ to $u$ by a unique
path of length $h$. Since they are in different branches of $T'$ they are connected in $T'$ by a path of length $2h$.
Together with $P$, this path forms a cycle of length $2h + 2r$ in $G$.  This works for any $r \in \{1,2,\dots,k\}$, giving cycles $C_{2h + 2},C_{2h+2},\dots,C_{2h+2k}$ and completing the proof. \qed

\section{Concluding remarks}

$\mbox{ }$ \vspace{-0.2in}

$\bullet$ The case of cycles of length zero mod $k$ has received considerable attention. Using his subdivided grid theorem for graphs of large tree-width, Thomassen~\cite{Thomassen} gave a polynomial-time algorithm for finding a cycle of length 0 mod $k$ in a graph or a certificate that no such cycle exists. It is an open question as to whether such an algorithm exists for finding a cycle of length $\ell$ mod $k$ when $\ell \neq 0$.

\medskip

$\bullet$ Thomassen~\cite{Thomassen} conjectured
that every graph of minimum degree at least $k + 1$ contains cycles of all possible even lengths mod $k$. This conjecture was proved
when $k$ is even by Liu and Ma~\cite{LiuMa}, and they showed further that there are cycles of $k/2$ consecutive even lengths in this case.
Liu and Ma also showed that if $k$ is odd and $G$ is a graph of minimum degree at least $k + 5$, then $G$ contains cycles of all even lengths mod $k$. It is natural to conjecture the following strengthening of Thomassen's conjecture~\cite{Thomassen}:

\begin{conjecture}
If $t \geq 1$, and $G$ is a graph with a maximum number of edges that does not contain cycles of $t$ consecutive even lengths,
then every block of $G$ is a complete graph of order at most $2t + 1$.
\end{conjecture}

It is an exercise to see that the conjecture is true for $t = 1$,
and in fact the blocks in every extremal graph with an odd number of vertices are triangles.
In particular, the conjecture (with $t=k/2$) implies $c_k(\ell) \leq k + 1$ whenever $k$ and $\ell$ are even, and $c_k(2) = k + 1$.
The result of Liu and Ma~\cite{LiuMa} shows $c_k(\ell) \leq 2k$ when both $k$ and $\ell$ are even, since every graph
of average degree at least $2k$ contains a subgraph of minimum degree at least $k + 1$.
We refer the reader to~\cite{Vsurvey} for a survey of this and related extremal problems for cycles in graphs.

\medskip

$\bullet$ In this paper we showed $c_k(\ell) = O(d_k(\ell)) = O(\ell k^{1/\ell})$ when $\ell$ is even.
We conjectured for each fixed $\ell$ that $c_k(\ell) \sim d_k(\ell)$ as $k \rightarrow \infty$ (see Conjecture \ref{cd}).
It may even be true that for infinitely many $n$, the extremal $n$-vertex graph with no cycle of length 4 mod $k$
is a connected graph whose blocks are all extremal $C_4$-free graphs with $k + 3$ vertices. We propose the more
tractable problem of determine the asymptotic value of $c_k(4)$. In this case our conjecture states $c_k(4) \sim k^{1/2}$,
since $d_k(4) \sim k^{1/2}$ (see~\cite{Vsurvey}).

\medskip

$\bullet$ Concerning chromatic number and consecutive cycle lengths, it is an open question to determine for $k \geq 2$ the largest chromatic number $\chi_k$ of a graph which does not contain cycles of $k$ consecutive lengths. Lemma 9 in~\cite{KSV} shows that a graph of chromatic
number at least $4k$ contains cycles of $k$ consecutive lengths, so we deduce $\chi_k \leq 4k - 1$. With some additional
effort, one can show $\chi_k \leq 2k + 2$ for $k \geq 2$. On the other hand, a graph $G^*_k$ whose blocks are all cliques of order $k + 1$ does not contain cycles of $k$ consecutive lengths, so $k + 1 \leq \chi_k \leq 2k + 2$ for $k \geq 2$. It seems plausible that $G^*_k$ is the extremal construction for $k \geq 3$, and perhaps $\chi_k = k + 1$ for all $k \geq 2$:

 \begin{conjecture}
 For all $k \geq 2$, $\chi_k = k + 1$.
 \end{conjecture}

 Gy\'{a}rf\'{a}s~\cite{G} proved that for $k \geq 2$, if a graph $G$ does not contain cycles of $k$ distinct odd lengths, then it has chromatic number at most $2k - 1$ with equality for $k \geq 3$ if and only if $G = G^*_{2k-1}$. Perhaps the same example is extremal for cycles of $k$ consecutive odd lengths.

\begin{center}
{\bf Acknowledgement}
\end{center}

This work was carried out when the second author visited the Institute for Mathematical Research (FIM) of ETH Z\"urich. He would like to thank FIM for the hospitality and for creating a stimulating research environment.

\end{document}